\documentclass[11pt]{amsart}
\usepackage{tikz}
\usepackage{booktabs}                                 
\usepackage[english]{babel}
\usepackage{latexsym}
 \usepackage[utf8]{inputenc}
\usepackage{babel}
\usepackage{fancyhdr}
\usepackage{longtable}
\usepackage{mathrsfs}
\usepackage{geometry}
\usepackage{textcomp}
\usepackage{caption}
\usepackage{lmodern}
\usepackage{mathrsfs}
 \usepackage[T1]{fontenc}
\usepackage[all,cmtip]{xy}
\usepackage{epsfig}
 \usepackage{amsmath,amsfonts,amssymb,amsthm}
\usepackage{graphics}
\usepackage{graphicx}
\usepackage{verbatim}
\usepackage{dsfont}
\usepackage{enumerate}  

\usepackage{pdflscape}
\usepackage{longtable}
\usepackage{booktabs}
\usepackage{colortbl}
\usepackage[shortlabels]{enumitem}

\newcommand{\C}{\mathbb{C}}

\newcommand{\QQ}{\mathbb{Q}}

\newcommand{\PP}{\mathbb{P}}

\newtheorem{thm}{Theorem}[section]

\newtheorem{rmk}[thm]{Remark}

\newtheorem*{aim*}{Aim of this paper}

\makeatletter    
\def\l@subsection{\@tocline{1}{0,2pt}{2pc}{8mm}{\ \ }} 
\makeatletter    
\def\l@section{\@tocline{1}{0,2pt}{2pc}{8mm}{\ \ }} 

\author{Lev Borisov}
\address{Department of Mathematics\\ Rutgers University\\ Piscataway, NJ 08854, USA}
\email[L.~Borisov]{borisov@math.rutgers.edu}
\author{Anders Buch}
\address{Department of Mathematics\\ Rutgers University\\ Piscataway, NJ 08854, USA}
\email[A.~Buch]{asbuch@math.rutgers.edu}

\author{Enrico Fatighenti }
\address{Department of Mathematical Sciences, Loughborough University, LE113TU, UK}
\email[E.~Fatighenti]{e.fatighenti@lboro.ac.uk}

\title[A journey from the octonionic $\mathbb P^2$ to a fake $\mathbb P^2$]{A journey from the octonionic $\mathbb P^2$ to a fake $\mathbb P^2$}

\begin{document}
\begin{abstract}
We discover a family of surfaces of general type with $K^2=3$ and $p=q=0$ as free $C_{13}$ quotients of special linear cuts of the octonionic projective plane $\mathbb O \mathbb P^2$. A special member of the family has $3$ singularities of type $A_2$, and is a quotient of a fake projective plane. We use the techniques of \cite{BF20} to define this fake projective plane by explicit equations in its bicanonical embedding.
\end{abstract}

\maketitle

\section{Introduction}
Fake projective planes are complex projective surfaces of general type with Hodge numbers equal to those of the usual projective plane $\C\PP^2$.
There are exactly $50$ complex conjugate pairs, constructed as ball quotients in \cite{CS11}
and they are fascinating gemstones in the vast mine of algebraic surfaces of general type.  The first explicit equations of a pair of fake projective planes were constructed in \cite{BK19}, and additional six pairs were given explicitly in \cite{BF20}. We refer the reader to \cite{BK19} for more background and history.

\smallskip
Many fake projective planes $\PP^2_{fake}$ admit an action of the cyclic group $C_3$. The quotient $\PP^2_{fake}/C_3$ is then a singular surface with $K^2=3$ and three singular points of type $A_2$. It can be deformed to construct interesting smooth surfaces with 
 $K^2=3$, genus $p=0$, and irregularity $q=0$. In \cite{BF20} the process was reversed and since the current paper is in many ways analogous, we describe \cite{BF20}  in some detail below.

\smallskip
The paper \cite{BF20} first builds a family of special complete intersections of seven Pl\"ucker hyperplanes in the Grassmannian $\operatorname{Gr}(3,\C^6)$ which admit a free action of the cyclic group $C_{14}$. This gives a family of surfaces $W$ which 
has $K_W^2=3, ~p=q=0$.  Then the authors find an element of this family such that the quotient by $C_{14}$ has an additional $C_3$ symmetry and three $A_2$ singularities. Its Galois cover (that was not at all easy to construct) is a fake projective plane with the automorphism group $(C_3)^2$, 
labeled by $(\operatorname{C2},p=2,\emptyset ,d_3D_3)$  in Cartwright-Steger classification \cite{CS11+}.

\smallskip
In the table \cite[Table 1]{BCP11} there are listed surfaces with fundamental group $C_{13}$ instead of $C_{14}$. The current paper is the result of our efforts to replicate the approach of \cite{BF20} and to construct their fundamental covers as complete intersections in some homogeneous space. We were not quite able to do it, instead we constructed them as \emph{almost} complete intersections of the $16$ dimensional octonionic projective plane $\mathbb O\PP^2$ in $\PP^{26}$ 
by certain $15$ linear equations, equivariant with respect to an order $13$ element in the Cartan subgroup of the $E_6$ group of automorphisms of $\mathbb O\PP^2$. 

\smallskip
Afterwards, the process was rather similar to that of \cite{BF20}, although there were some technical complications due to lack of unramified double covers. In particular, we had more difficulty controlling the size of the coefficients and had to work with $60K$ decimal digit numbers at some intermediate steps.

\smallskip
The paper is organized as follows.
In Section 2 we describe our motivation for using the octonionic projective plane and the special linear cuts that achieve our goal of constructing surfaces with $K^2=3,~p=q=0$. We also describe how we found a special element of this family with three $A_2$ singularities.  In Section 3 we briefly explain the construction of the fake projective plane $\PP^2_{fake}$, labeled by $(\operatorname{C18},p=3,\emptyset,d_3 D_3)$ in \cite{CS11+}, and state some open problems.

\smallskip
{\bf Acknowledgements.} Our computations relied heavily on Mathematica software package \cite{Ma} and to some extent on Julia \cite{Julia},  Macaulay2 \cite{M2} and Magma \cite{Mag}. We thank John Cremona for helping us access the number theory server at the University of Warwick.
Anders Buch was partially supported by the NSF grant DMS-1503662.

\section{Special cuts of the octonionic projective plane}
\subsection{Motivation.} As was mentioned in the Introduction, we set out to find a family of surfaces of general type with $K^2=3, ~p=q=0$ and fundamental group $C_{13}$. Here is how this search led us to consider cuts of the octonionic projective plane $\mathbb O\PP^2$.

\smallskip
Let $X$ be the universal cover of a surface in question, with a free action of an order $13$ automorphism $g$. Then $K_X^2=39$ and $\chi(K_X)=13$. It is reasonable to expect that $h^1(X,K_X)=0$ and $h^0(X,K_X)={12}$. Then the  pluricanonical ring $\bigoplus_{n\geq 0} H^0(X,nK_X)$ of $X$ must have the  graded dimension
$$
\sum_{n\geq 0} \dim H^0(X,nK_X) \,t^n =  1 + 12 t + 52 t^2 + 130 t^3 + \cdots = \frac {1 + 9 t + 19 t^2 + 9 t^3 + t^4}{(1-t)^3}.
$$
It is also reasonable to assume that the pluricanonical ring is generated in degree one, so $X$ is embedded into $\PP^{11}$. It is also plausible that its image is cut out by $12 (12+1)/2 - 52 = 26$ quadrics.
By the Holomorphic Lefschetz formula, as in \cite[Theorem 2.1]{Hua11}, the trace of the action of $g$ on $H^0(X,K_X)$ is $(-1)$, which  means that $H^0(X,K_X)$  has a basis of eigenvectors of $g$ with eigenvalues $\zeta_{13}^i$ for $i=1,\ldots, 12$. Similarly, the action of $g$ on the space of quadrics splits it into $13$ two-dimensional eigenspaces.

\smallskip
Inspired by \cite{BF20}, we undertook a rather exhaustive computer search for homogeneous varieties of degree $39$ and other relevant invariants, but were not successful. However, the octonionic projective plane  $\mathbb O\PP^2$  has degree $78$, and we observed the following remarkable coincidence: the homogeneous coordinate ring of  $\mathbb O\PP^2$ has graded dimension
$$
\frac {(1 + 9 t + 19 t^2 + 9 t^3 + t^4)(1+t)}{(1-t)^{17}}.
$$
More specifically, $\mathbb O\PP^2$ is the dim $16$ singular locus of the $E_6$-invariant cubic in $\PP^{26}$ cut out by the $27$ quadratic equations which are the partial derivatives of the cubic. So our idea was to take a linear cut of $\mathbb O\PP^2$ by $15$ equations (so that we are in $\mathbb P^{11}$) which only drop the dimension by $14$. We also want one of the quadratic equations to reduce to zero on the linear cut.
The best analogy would be cutting a quadric $xy=zw$ with Hilbert series $\frac {1+t}{(1-t)^3}$ by two linear equations $x=0$ and $z=0$ to get a line with the Hilbert series $\frac 1{(1-t)^2}$ but, ultimately, it was a lucky guess.

\smallskip
\subsection{Octonionic projective plane.}
There are several incarnations of the $E_6$-invariant cubic found in the literature. We used the one in Jacob Lurie's undergraduate thesis \cite{Lurie}, namely 
$$
\begin{array}{l}
-P_{10} P_{13} P_{16} - P_{11} P_{14} P_{17} - P_{12} P_{15} P_{18} - P_{16} P_{17} P_{18} + P_{1} P_{10} P_{19} - 
 P_{1} P_{18} P_{2} + P_{11} P_{2} P_{20} - P_{1} P_{14} P_{21} 
 \\
 - P_{16} P_{20} P_{21}
  - P_{18} P_{22} P_{23} - 
 P_{17} P_{19} P_{24} - P_{13} P_{2} P_{24} - P_{14} P_{15} P_{25} - P_{19} P_{22} P_{25} - P_{12} P_{13} P_{26} 
 - 
 P_{20} P_{23} P_{26}
 \\ - P_{10} P_{11} P_{27} - P_{21} P_{24} P_{27} - P_{25} P_{26} P_{27} + P_{12} P_{21} P_{3} - 
 P_{10} P_{23} P_{3} - P_{2} P_{25} P_{3} - P_{15} P_{20} P_{4} + P_{13} P_{22} P_{4} 
 \\- P_{17} P_{3} P_{4} - 
 P_{12} P_{19} P_{5} + P_{14} P_{23} P_{5} - P_{27} P_{4} P_{5} - P_{11} P_{22} P_{6} + P_{15} P_{24} P_{6} - 
 P_{1} P_{26} P_{6} - P_{16} P_{5} P_{6} 
 \\- P_{11} P_{12 }P_{7} 
 - P_{23} P_{24} P_{7} + P_{16} P_{25} P_{7} - 
 P_{1} P_{4} P_{7} - P_{10} P_{15} P_{8} - P_{21} P_{22 }P_{8}
  + P_{17} P_{26} P_{8} - P_{2} P_{5} P_{8}
  \\ - 
 P_{13} P_{14} P_{9} - P_{19} P_{20} P_{9} + P_{18} P_{27} P_{9} - P_{3} P_{6} P_{9 }- P_{7} P_{8 }P_{9}
\end{array}
$$
The $27$ variables $P_1,\ldots, P_{27}$ are indexed by the lines on the Fermat cubic surface in $\C\PP^3$ and the terms correspond to triples of coplanar lines. The sign prescription is more intricate, given in terms of the $C_3$ action on the cubic,
see \cite{Lurie}. The octonionic projective plane $\mathbb O\PP^2$ is cut out by the $27$ partial derivatives of the above cubic.

\smallskip
There is a  Cartan subgroup $(\C^*)^6$ of $E_6$ that acts diagonally on the variables $P_i$. We picked an element $g$ of order $13$ of it which acts by  $P_i\mapsto \zeta_{13}^{a_i}P_i$ with the weights $a_i $ given by 
$$
(6, 7, 7, 10, 3, 6, 10, 3, 0, 5, 8, 8, 4, 9, 5, 4, 9, 0, 2, 11, 11, 
12, 1, 2, 12, 1, 0).
$$
As the reader can see, the action of $g$ on the variables has a three-dimensional eigenspace of weight zero and $12$ two-dimensional eigenspace of other weights.
For the three invariant variables $P_9,P_{18},P_{27}$, the corresponding partial derivatives of the cubic 
\begin{equation}\label{invquad}
\begin{array}{l}
-P_{13} P_{14} - P_{19} P_{20} + P_{18} P_{27 }- P_{3} P_{6 }- P_{7} P_{8}, 
\\ -P_{12} P_{15} - P_{16} P_{17} - 
  P_{1} P_{2 }- P_{22 }P_{23} + P_{27} P_{9}, 
  \\ -P_{10} P_{11} - P_{21} P_{24} - P_{25} P_{26} - P_{4} P_{5} + 
  P_{18} P_{9}
\end{array}
\end{equation}
involve all of the variables $P_i$.

\smallskip
At this point, our expectations of the $g$-action on  $H^0(X,K_X)$  indicate that we need to take a linear cut by
\begin{align*}
(P_{9}, P_{18}, P_{27}, P_{23} + d_1 P_{26}, P_{19} + d_2 P_{24}, P_{5} + d_3 P_{8}, 
   P_{13} + d_4 P_{16}, P_{10} + d_5 P_{15}, P_{1} + d_6 P_{6},
   \\
    P_{2} + d_7 P_{3}, P_{11} + d_8 P_{12}, 
   P_{14} + d_9 P_{17}, P_{4} + d_{10} P_{7}, P_{20} + d_{11} P_{21}, P_{22 }+ d_{12} P_{25})
\end{align*}
for some constants $d_1,\ldots, d_{12}$.
Moreover, we want a linear combination of the $g$-invariant quadrics  \eqref{invquad} to vanish on the linear subspace of the cut.
In view of the Cartan subgroup symmetry, it is reasonable to pick this linear combination to be the sum of the above quadrics. This gives $6$ simple equations on $d_i$, namely $d_i d_{13-i}= -1$. We have been able to verify by computer at a specific point that the resulting scheme is a smooth surface of degree $13$ and is thus a good candidate for our $X$.

\smallskip
Specifically, the $26$  quadrics that cut out  $X$ are given by 
\begin{equation}\label{bad26}
\begin{array}{l}
 -t_{10}^2 + d_2 t_2 t_5 - t_1 t_6 - d_2 t_{11} t_9,\  
  t_3^2 + t_2 t_4 - d_1 t_{12} t_7 + t_{11} t_8,\  
  -d_1 t_1 t_5 - d_1 t_{12} t_7 + d_2 t_{11} t_8 - t_{10} t_9,\  
  \\
  -t_{12} t_4 - t_{11} t_5 + t_{10} t_6 - t_7 t_9,\  
  -t_4 t_6 - t_3 t_7 + d_2 t_2 t_8 - d_1 t_1 t_9,\  
  t_3 t_4 + t_2 t_5 + t_1 t_6 - t_{12} t_8,\  
  \\
  d_1 t_1 t_2 + d_1 t_{12} t_4 + t_{10} t_6 - t_8^2,\  
  -d_2 t_{11} t_{12} + t_5^2 + t_3 t_7 + t_1 t_9,\  
  d_2 t_{11} t_2 - t_{10} t_3 - t_6 t_7 + t_4 t_9,\  
  \\
  t_4^2 - t_3 t_5 + d_2 t_2 t_6 + d_1 t_1 t_7,\  
  -t_{12} t_6 + t_{11} t_7 - t_{10} t_8 - t_9^2,\  
  -d_2 t_2 t_3 + t_1 t_4 + d_2 t_{11} t_7 - t_{10} t_8,\  
  \\
  t_{10} t_{12} + t_4 t_5 - t_2 t_7 - t_1 t_8,\  
  d_1 t_1 t_3 - d_1 t_{12} t_5 + d_2 t_{11} t_6 - t_8 t_9,\  
  -t_{10} t_{11} + t_3 t_5 + t_2 t_6 - d_1 t_{12} t_9,\  
  \\
  -d_2 t_{11}^2 + d_1 t_{10} t_{12} - t_4 t_5 + t_3 t_6,\  
  d_2 t_2^2 + t_1 t_3 - t_{10} t_7 - t_8 t_9,\  
  d_1 t_1 t_{12} + t_6 t_7 - t_5 t_8 - t_4 t_9,\  
  \\
  -d_1 t_{12}^2 + t_5 t_6 + t_3 t_8 - t_2 t_9,\  
  d_1 t_1^2 - d_2 t_{11} t_4 - t_{10} t_5 + t_7 t_8,\  
  -t_{12} t_3 - t_{11} t_4 + t_7 t_8 + t_6 t_9,\  
  \\
  d_1 d_2 t_{12} t_2 - d_2 t_{11} t_3 - t_{10} t_4 - t_6 t_8,\  
  -t_1 t_{11} - t_{10} t_2 + t_5 t_7 - t_3 t_9,\  
  d_1 t_1 t_{10} + t_5 t_6 + t_4 t_7 + d_2 t_2 t_9,\  
  \\
  d_2 t_{12} t_2 + t_{10} t_4 - t_7^2 - t_5 t_9,\  
  d_1 t_1 t_{11} + t_6^2 + t_4 t_8 + t_3 t_9
  \end{array}
\end{equation}
in the homogeneous coordinates $(t_1:\ldots:t_{12})$ of $\PP^{11}$. The action of $g$ is $t_i\mapsto \zeta_{13}^i t_i$.

\begin{rmk}
The action of the Cartan subgroup of $E_6$ reduces the dimension of the space of parameters $d$ from six to two (taking into account the need to preserve the invariant quadric that has to vanish on the cut reduces $E_6$ to $F_4$). We expect the total family to have dimension four, but it is not clear how one can build it. What makes the elements above special is that these surfaces $X$ admit an additional $C_3$ symmetry that extends the $C_{13}$ action to the semidirect product of these two groups. Namely, by scaling the variables (but still calling them $t_i$) we could rewrite the equations  \eqref{bad26} as 
\begin{equation}\label{good26}
\begin{array}{l}
-t_{10}^2 - d_1 d_2^2 (t_2 t_5 + t_1 t_6 - t_{11} t_9),\  
  d_1 d_2^2 t_3^2 + t_2 t_4 + t_{12} t_7 - t_{11} t_8,\  
  d_1 d_2^2 t_1 t_5 + t_{12} t_7 - d_2 (t_{11} t_8 + t_{10} t_9),\ 
  \\ -t_{12} t_4 
  + 
   d_1 d_2 (-t_{11} t_5 + t_{10} t_6 + d_2 t_7 t_9),\  
  t_4 t_6 + d_2 (-t_3 t_7 - t_2 t_8 + d_1 d_2 t_1 t_9),\  
  d_1 d_2 (t_3 t_4 - t_2 t_5 + d_2 t_1 t_6) - t_{12} t_8,\  
  \\
  d_1 d_2^2 t_1 t_2 + t_{12} t_4 + d_2 t_{10} t_6 - d_2 t_8^2,\  
  t_{11} t_{12} + d_1 d_2 (t_5^2 - t_3 t_7 + d_2 t_1 t_9),\  -t_{11} t_2 - t_{10} t_3 + 
   t_6 t_7 + t_4 t_9,\  
   \\
  t_4^2 + d_1 d_2^2 (t_3 t_5 + t_2 t_6 - t_1 t_7),\  -t_{12} t_6 + t_{11} t_7 - t_{10} t_8 - 
   d_1 d_2^2 t_9^2,\  -d_1 d_2^2 t_2 t_3 + d_2 t_1 t_4 + d_2 t_{11} t_7 - t_{10} t_8,\  
   \\
  t_{10} t_{12} - d_1 d_2 (t_4 t_5 - t_2 t_7 + d_2 t_1 t_8),\  
  d_1 d_2^2 t_1 t_3 + t_{12} t_5 - d_2 (t_{11} t_6 + t_8 t_9),\  
  t_{10} t_{11} - d_1 d_2 (d_2 t_3 t_5 - t_2 t_6 + t_{12} t_9),\ 
  \\ -d_2 t_{11}^2 + t_{10} t_{12} + 
   d_2 t_4 t_5 + d_1 d_2^2 t_3 t_6,\  t_{10} t_7 + d_1 d_2 (t_2^2 + d_2 t_1 t_3 - t_8 t_9),\ 
   t_1 t_{12} - t_6 t_7 + t_5 t_8 - t_4 t_9,\  
   \\
   -t_{12}^2 - 
   d_1 d_2^2 (t_5 t_6 - t_3 t_8 + t_2 t_9),\  
  d_1 d_2^2 t_1^2 + t_{11} t_4 + t_{10} t_5 - t_7 t_8,\  
  t_{11} t_4 - d_2 (t_{12} t_3 + t_7 t_8) + d_1 d_2^2 t_6 t_9,\ 
  \\ -t_{10} t_4 + 
   d_1 d_2 (t_{12} t_2 + d_2 t_{11} t_3 - t_6 t_8),\  
  d_2 t_1 t_{11} - t_{10} t_2 + d_2 t_5 t_7 - d_1 d_2^2 t_3 t_9,\  -t_4 t_7 + 
   d_1 d_2 (t_1 t_{10} - t_5 t_6 + d_2 t_2 t_9),\  
   \\
  d_2 t_{12} t_2 + t_{10} t_4 - d_2 t_7^2 + d_1 d_2^2 t_5 t_9,\  
  t_4 t_8 + d_1 d_2 (-t_1 t_{11} + t_6^2 + d_2 t_3 t_9)
\end{array}
\end{equation}
with the additional symmetry $t_{i}\mapsto t_{3 i \hskip -4pt \mod 13}$.
The details are in \cite[Section2.nb]{BBF20+}.
\end{rmk}

\subsection{Constructing a cut with $A_2$ singularities.}
Our method of constructing a fake projective plane largely followed the blueprint of \cite{BF20}. 

\smallskip
We set $d_3=d_4=d_5=d_6=1$ and tried to find out which $(d_1,d_2)$ give singular cuts. In order to achieve this, we worked on an affine coordinate chart of $\mathbb O\PP^2$ which can be obtained by solving the equations of $\mathbb O\PP^2$ 
for eleven of the variables as follows.
$$
\begin{array}{l}
 P_{4} =P_{10} P_{16} + P_{2} P_{24} + P_{12} P_{26} + P_{14} P_{9}, 
 P_{6} =-P_{14} P_{17} + P_{2} P_{20} - P_{10} P_{27} - P_{12} P_{7}, 
\\
 P_{8} =-P_{1} P_{14} - P_{16 }P_{20 }- P_{24} P_{27} + P_{12} P_{3},
P_{11} = P_{15} P_{24} - P_{1} P_{26 }- P_{16} P_{5 }- P_{3} P_{9},
\\
 P_{13} =P_{15 }P_{20} + P_{17} P_{3 }+ P_{27} P_{5 }+ P_{1} P_{7}, 
 P_{18} =-P_{20} P_{26} - P_{10 }P_{3} + P_{14 }P_{5} - P_{24 }P_{7}, 
 \\
 P_{19} =-P_{14} P_{15} - P_{26 }P_{27 }- P_{2} P_{3} + P_{16} P_{7}, 
 P_{21} =-P_{10} P_{15} + P_{17 }P_{26} - P_{2} P_{5} - P_{7} P_{9}, 
  \\
  P_{22}=1, 
  P_{23} =-P_{12} P_{15} - P_{16} P_{17} - P_{1} P_{2} + P_{27} P_{9}, 
 P_{25} =P_{1} P_{10} - P_{17} P_{24} - P_{12 }P_{5 }- P_{20 }P_{9}
 \end{array}
$$
We obtained this chart by connecting the formulas for the Cartan cubic from \cite{GrossElkies} and \cite{Lurie}.
We then further solved for five of the variables to reduce their number while still keeping the equations relatively short.
Then we looked for tangent vectors for the surfaces with $d_1=1$ that lie in a codimension three subspace, by a multivariable
Newton method starting at random points. The idea is that some of these would happen at values of $d_2$
where the surface $X=X_{1,d_2}$ acquires a node. After some trial and error we saw that solutions to
$$
-27 - 34 d_2 - 397 d_2^2 - 172 d_2^3 - 821 d_2^4 + 190 d_2^5 - 83 d_2^6 + 16 d_2^7 = 0
$$
give singular surfaces. As in \cite{BF20}, we then perturbed $d_1$ slightly to $1 + 10^{-20}$ to find a nearby point on the locus of singular surfaces. This lead us to conjecture that generic points $(d_1,d_2)$ on the curve
$$
\begin{array}{l}
0=-4 d_1^3 + 8 d_1^4 - 4 d_1^5 - 12 d_1^2 d_2 - 16 d_1^3 d_2 
+ 28 d_1^4 d_2 - 
 39 d_1^5 d_2 + 12 d_1^6 d_2 - 12 d_1 d_2^2 - 28 d_1^2 d_2^2 - 54 d_1^3 d_2^2 
 \\+ 
 78 d_1^4 d_2^2 - 34 d_1^5 d_2^2 
 + 28 d_1^6 d_2^2 - 12 d_1^7 d_2^2 - 4 d_2^3 - 
 39 d_1 d_2^3 - 34 d_1^2 d_2^3 - 277 d_1^3 d_2^3 + 192 d_1^4 d_2^3 - 
 277 d_1^5 d_2^3
 \\ + 54 d_1^6 d_2^3 - 16 d_1^7 d_2^3 + 4 d_1^8 d_2^3
  - 8 d_2^4 - 
 28 d_1 d_2^4 - 78 d_1^2 d_2^4 - 192 d_1^3 d_2^4 + 192 d_1^5 d_2^4 - 
 78 d_1^6 d_2^4 + 28 d_1^7 d_2^4 
 \\- 8 d_1^8 d_2^4 - 4 d_2^5 - 16 d_1 d_2^5 - 
 54 d_1^2 d_2^5 - 277 d_1^3 d_2^5 - 192 d_1^4 d_2^5 - 277 d_1^5 d_2^5 + 
 34 d_1^6 d_2^5 - 39 d_1^7 d_2^5 + 4 d_1^8 d_2^5 + 12 d_1 d_2^6 
 \\+ 
 28 d_1^2 d_2^6 + 34 d_1^3 d_2^6 + 78 d_1^4 d_2^6 + 54 d_1^5 d_2^6 - 
 28 d_1^6 d_2^6 + 12 d_1^7 d_2^6 - 12 d_1^2 d_2^7 - 39 d_1^3 d_2^7 - 
 28 d_1^4 d_2^7 - 16 d_1^5 d_2^7 
 \\+ 12 d_1^6 d_2^7 + 4 d_1^3 d_2^8 + 
 8 d_1^4 d_2^8 + 4 d_1^5 d_2^8
 \end{array}
$$
give nodal $X_{d_1,d_2}$. 

\smallskip
We then looked for singular points of this curve. There were several such points, one of which was a cusp of the curve. We focused our attention on it and discovered a surface $X_{d_1,d_2}$ with $39$ $A_2$ singularities. Specifically, both $d_1$ and $d_2$ can be given given as roots of 
$$
0=2187 + 7290 d + 23433 d^2 + 21640 d^3 + 66393 d^4 - 21640 d^5 + 23433 d^6 - 7290 d^7 + 2187 d^8 
$$
approximately given by $(d_1,d_2)\approx (1.93 + 2.30\, \rm i, 0.0125 - 0.515\,\rm i)$. Of course, the same is true for all of the Galois conjugates of this pair. From now on we will call this surface $X_0$.

\smallskip
We observed that four of the Galois conjugate pairs of $(d_1,d_2)$ give isomorphic surfaces. To see that, we noticed that scheme 
$X_{d_1,d_2}$ cut out by \eqref{good26} is isomorphic to $X_{-1/d_2,d_1}$ under the coordinate change 
\begin{equation}\label{4gal}
(t_1,\ldots, t_{12}) \mapsto ( d_2 t_5, t_{10},  d_2 t_2,  -d_1 t_7, 
  t_{12}, t_4,-d_1 d_2 t_9,  -d_1 d_2 t_1,  d_2 t_6,
    -d_1 t_{11},  -d_1 d_2 t_3,  -d_1 t_8)
\end{equation}
The idea behind it was to use $i\to 5 i \hskip -4pt \mod 13$, and we heavily relied on Mathematica computations, see \cite[Section2.nb]{BBF20+}.

\smallskip
We then used the symmetry \eqref{4gal} to average the $C_{13}$-invariants of the coordinate ring of the surface $X_0$ suspected to have $A_2$ singularities, to get $X_0/C_{13}$ defined in the $13$-dimensional weighted projective space $W\PP(2^4,3^{10})$ by $9$ equations of degree five  and $29$ equations of degree six. Due to the above symmetrization, the coefficients were in the field $\QQ(\sqrt {-2})$. 

\subsection{Finding singular points.} It was not entirely trivial to find the singular points of $X_0/C_{13}$. We did it by calculating a degree $12$ equation in the first four variables which gives a (non-normal) image of $X_0/C_{13}$ in $\PP^3$. Then we looked for its curves of singularities by finding multiple singular points on  random hyperplane cuts. Then we have looked for singular points outside of the curve of singularities, and indeed hit upon $A_2$ singularities. We were then able to verify that these were the only singularities by computing the degree of the singular locus over a finite field. As in the case of \cite{BF20}, the $A_2$ singularities were not defined over the quadratic extension of $\QQ$, but a coordinate change gave us a model of $X_0/C_{13}\subseteq W\PP^{13}$ still defined over  $\QQ(\sqrt{-2})$ and with three singular points defined over $\QQ$.

\smallskip

\section{Constructing the fake projective plane}
\subsection{Constructing the triple cover.}
By the work of Keum \cite{Keum.comment}, the surface $X_0/C_{13}$ admits a Galois triple cover which is a fake projective plane. In this, it is very similar to the situation in \cite{BF20} and we employed the same general method. It was useful that in both cases there was an additional order three automorphism $\sigma$  because the FPP had a $C_3\times C_3$ group of automorphisms. Specifically, we looked for sections $f$ and $d$ of  $4K_{X_0/C_{13}}$ which satisfy
$$
f \,\sigma(f) \sigma^2(f) = d^3.
$$
Moreover $f$ and $d$ should have certain vanishing on the exceptional lines at the blowup of $A_2$ singularities. We refer the reader to \cite{BF20} for details.

\smallskip
The nature of $X_0/C_{13}$ made the computations more challenging. In particular, at some point we had to work with random points on the surface computed with $6\times 10^4$ digits of accuracy. The equations for $f$ and $d$ had coefficients in $\QQ(\sqrt{-2})$ which were about $1.5 \times 10^4$ digits long. As in \cite{BF20}, we solved it over a finite field of $19$ elements, but now we used a p-adic version of the Newton's method to quickly gain the needed accuracy.

\smallskip
Once the triple cover was constructed, we used the fixed points of the automorphisms of $\PP^2_{fake}$ to get a basis with nicer equations, only about 100 digits long coefficients, see \cite[Section3.nb]{BBF20+} for details.
This surface is  labeled by $(\operatorname{C18},p=3,\emptyset,d_3 D_3)$ in the classification of \cite{CS11+}, since it is the only one with an automorphism group that contains $(C_3)^2$ and Picard group that contains $C_{13}$.

\smallskip
The details of the above process are in \cite[Section3.nb]{BBF20+}.

\smallskip
\subsection{Open questions.} Let us now discuss open problems related to this construction.

\smallskip
The first question is how to verify that the special cuts $X_{d_1,d_2}$  of $\mathbb O\PP^2$
are simply connected. Since these are not complete intersections, the Lefschetz Hyperplane theorem can not be applied, so other methods are needed. It might perhaps follow from our construction and \cite{CS11+}, but a more direct argument is desirable.

\smallskip
A related question is how to construct non-$C_3$-invariant deformations of $X_{d_1,d_2}$. It looks like they will no longer be cuts of $\mathbb O\PP^2$ but perhaps one can get them by carefully examining the equations  \eqref{bad26}. 

\smallskip
The quotient of the fake projective plane $(\operatorname{C18},p=3,\emptyset,d_3 D_3)$ by $(C_3)^2$ is also covered by $(\operatorname{C18},p=3,\{2I\})$.
This fake projective plane in turn covers a surface $(\operatorname{C18},p=3,\{2\})$, which is covered by three other fake projective planes.
The method of \cite{BF20} is, unfortunately, not quite applicable here, so how do we find (the equations of) these other surfaces?

\smallskip
It is known from \cite{CS11+} that $(\operatorname{C18},p=3,\emptyset,d_3 D_3)$ has Picard group $C_2\times C_2\times C_{13}$. While the $C_{13}$ part can  be inferred from our construction (even though it may not be entirely trivial to follow), the other two factors are mysterious. It would be interesting to see them explicitly, and they may be useful in answering both the previous and the next questions.

\smallskip
A perennial question is how one can reduce the size of the coefficients in the equations. There are currently only ad hoc tools that are not very successful, except in \cite{BK19} case.

\end{document}